\documentclass[letterpaper, 10pt, conference]{ieeeconf}
\IEEEoverridecommandlockouts 
\usepackage{amsmath,amssymb,url}
\usepackage{graphicx,subfigure,hyperref}
\usepackage{color,comment}

\newcommand{\norm}[1]{\ensuremath{\left\| #1 \right\|}}

\newcommand{\braces}[1]{\ensuremath{\left\{ #1 \right\}}}
\newcommand{\parenth}[1]{\ensuremath{\left( #1 \right)}}

\newcommand{\refeqn}[1]{(\ref{eqn:#1})}

\newcommand{\deriv}[2]{\ensuremath{\frac{\partial #1}{\partial #2}}}

\newcommand{\SO}{\ensuremath{\mathsf{SO(3)}}}
\newcommand{\T}{\ensuremath{\mathsf{T}}}

\renewcommand{\Re}{\ensuremath{\mathbb{R}}}
\newcommand{\Sph}{\ensuremath{\mathsf{S}}}

\title{\LARGE \bf
Global Formulations of Lagrangian and Hamiltonian Mechanics\\ on Two-Spheres}

\author{Taeyoung Lee\authorrefmark{1}, Melvin Leok\authorrefmark{2}, and N. Harris McClamroch%
\thanks{Taeyoung Lee, Mechanical and Aerospace Engineering, George Washington University, Washington DC 20052 {\tt tylee@gwu.edu}}
\thanks{Melvin Leok, Mathematics, University of California at San Diego, La Jolla, CA 92093 {\tt mleok@math.ucsd.edu}}%
\thanks{N. Harris McClamroch, Aerospace Engineering, University of Michigan, Ann Arbor, MI 48109 {\tt
nhm@umich.edu}}%
\thanks{\textsuperscript{\footnotesize\ensuremath{*}}This research has been supported in part by NSF under the grants CMMI-1243000 (transferred from 1029551), CMMI-1335008, and CNS-1337722.}
\thanks{\textsuperscript{\footnotesize\ensuremath{\dagger}}This research has been supported in part by NSF under grants DMS-1010687, CMMI-1029445, DMS-1065972, CMMI-1334759, DMS-1411792, DMS-1345013.}
}

\newtheorem{prop}{Proposition}
\newtheorem{cor}{Corollary}

\begin{document}
\allowdisplaybreaks

\maketitle \thispagestyle{empty} \pagestyle{empty}

\begin{abstract}
This paper provides global formulations of Lagrangian and Hamiltonian variational dynamics evolving on the product of an arbitrary number of two-spheres. Four types of Euler-Lagrange equations and Hamilton's equations are developed in a coordinate-free fashion on two-spheres, without relying on local parameterizations that may lead to singularities and cumbersome equations of motion. The proposed intrinsic formulations of Lagrangian and Hamiltonian dynamics are novel in that they incorporate the geometry of two-spheres, resulting in equations of motion that are expressed compactly, and they are useful in analysis and computation of the global dynamics. 
\end{abstract}

\section{Introduction}

The two-sphere is the two-dimensional manifold that is composed of unit-vectors in $\Re^3$. There are a wide variety of dynamical systems that evolve on multiple copies of two-spheres. In robotics, the configuration of articulated robotic arms interconnected by spherical joints is represented by two-spheres~\cite{MurLi93}. Since the surface of the Earth is approximately a sphere, spherical dynamics arise readily in earth science and meteorology~\cite{New01}. 
In quantum mechanics, the pure state space of a two-level quantum mechanical system is a two-sphere, referred to as the Bloch sphere~\cite{BloPR46}. 

In most of the existing literature on dynamical systems evolving on two-spheres, the unit-sphere is parameterized by two angles. For example, a point on the two-sphere is often described by its longitude and latitude. Parametrizing the two-sphere by two angles is straightforward, and the angles are typically viewed as being in an open subset of $\Re^2$. 

However, such parameterizations of the two-sphere suffer from the following two main issues. First, parameterizations represent the two-sphere only locally. This can be easily observed from the fact that the longitude is not well defined at the north pole and the south pole. This causes a singularity in representing the kinematics on the two-sphere, especially when converting the velocity of a curve on the two-sphere into the time-derivatives of the longitude and the latitude. This yields numerical ill-conditioning in the vicinity of those singularities, which cannot be avoided unless one switches coordinate charts, which becomes problematic when trying to track motions with large angular deviations.


The second issue is that the equations of motion of dynamical systems on the two-sphere become exceedingly complicated when expressed using local coordinates and necessarily involve complicated trigonometric expressions.  For example, the dynamics of a multiple spherical pendulum, written in terms of angles, is extremely complicated. 

This paper aims to provide global formulations of dynamics evolving on the product of an arbitrary number of two-spheres. In particular, we study dynamical systems that can be viewed as Lagrangian systems or as Hamiltonian systems that encompass a large class of mechanical systems that appear in robotics, structural dynamics, quantum mechanics or meteorology. Most importantly, the unit-vectors on the two-sphere are regarded as elements of a manifold, and dynamics are formulated directly on two-spheres in a global fashion via variational principles.

This geometric formulation is said to be \textit{coordinate-free}, as it does not require the use of local charts\index{charts}, coordinates\index{coordinates} or parameters that may lead to singularities or ambiguities in the representation. As such, it can be applied to arbitrarily large angle rotational maneuvers on the two-spheres globally. Furthermore, this provides an efficient and elegant way to formulate,  analyze, and compute the dynamics and their temporal evolutions. The corresponding mathematical model developed on two-spheres is nicely structured and elegant. This representational efficiency has a substantial practical advantage compared with local coordinates for many complex dynamical systems; this fact has not been appreciated by the applied scientific and engineering communities.

In short, the main contribution of this paper is providing geometric formulations of the equations of motion for Lagrangian and Hamiltonian systems that evolve on two-spheres using variational methods. The proposed global formulations, that do not require local charts, have not been previously studied, even in the well-known literature on geometric mechanics, such as~\cite{Mar92,HolSch09,MarRat99}.   Preliminary results have been given in~\cite{LeeLeoIJNME09,LeeLeoPICDC12}, where Euler-Lagrange equations are developed for a certain class of mechanical systems whose kinetic energy is repressed as a quadratic form with fixed inertia elements. This paper provides both Euler-Lagrange equations and Hamilton's equations for arbitrary mechanical systems without such restrictions. 

\section{The Two-Sphere}

The two-sphere is the two-dimensional manifold of unit-vectors in $\Re^3$, i.e.,
\begin{align}
\Sph^2=\{q\in\Re^3\,|\, \|q\|=1\}.
\end{align}
It is composed of the set of points that have the unit distance from the origin in $\Re^3$. The tangent space of the two-sphere at $q\in\Sph^2$ corresponds to the two-dimensional plane that is tangent to the sphere at the point $q$, and it is given by
\begin{align}
\T_q\Sph^2 =\{\xi\in\Re^3\,|\, q\cdot \xi =0\}.
\end{align}
Throughout this paper, the standard dot product between two vectors is denoted by $x\cdot y = x^T y$ for any $x,y\in\Re^n$. 

Therefore, for any curve $q(t):\Re\rightarrow\Sph^2$ on the two-sphere parameterized by time $t$, its time derivative satisfies $q(t)\cdot\dot q(t)=0$. From now on, we do not explicitly denote dependence on time for brevity, unless needed. This implies that there exists an angular velocity $\omega:\Re\rightarrow\Re^3$ such that
\begin{align}
\dot q = \omega\times q = S(\omega) q,\label{eqn:rotkin5}
\end{align}
where the hat map $S(\cdot):\Re^3\rightarrow \Re^{3\times 3}$ is defined such that $S(x)y=x\times y$ and $S^T(x)=-S(x)$ for any $x,y\in\Re^3$. More explicitly, 
\begin{align}
S(\omega) = \begin{bmatrix} 0 & -\omega_3 & \omega_2 \\
\omega_3 & 0 & -\omega_1\\
-\omega_2 & \omega_1 & 0
\end{bmatrix}.
\end{align}
Without loss of generality, the angular velocity is constrained to be orthogonal to $q$, i.e., $\omega\cdot q=0$. Therefore, the three vectors $q$, $\dot q$, and $\omega$ are mutually orthogonal, and the angular velocity can be written as 
\begin{align}
\omega=S(q)\dot q. \label{eqn:w}
\end{align}
It follows that $\dot\omega = S(q)\ddot q$ and $\dot\omega$ is perpendicular to $q$ as well.

\section{Lagrangian Mechanics on Two-Spheres}

We consider dynamical systems evolving on the product of $n$ copies of two-spheres, namely $(\Sph^2)^n$. The corresponding Euler--Lagrange equations are first obtained in terms of a Lagrangian expressed in terms of the configurations and the time derivatives of the configurations, namely $(q,\dot q)$.   A second form of Euler--Lagrange equations is obtained in terms of a modified Lagrangian expressed in terms of the configurations and the angular velocities, $(q,\omega)$. In each case, these Euler--Lagrange equations are simplified for the important case that the kinetic energy function is a quadratic function of the time derivatives of the configurations.  

\subsection{Euler--Lagrange equations in terms of $(q,\dot q)$}

Suppose that a Lagrangian $L(q,\dot q):\T(\Sph^2)^n\rightarrow\Re^1$ is given on the tangent bundle of $(\Sph^2)^n$, where $(q,\dot q) = ((q_1,\ldots, q_n),(\dot q_1,\ldots, \dot q_n))\in \T(\Sph^2)^n$. For example, it can be defined as the difference between a kinetic energy and a potential energy. We derive the corresponding Euler--Lagrange equations according to Hamilton's variational principle.

Let $q_i: [t_0,t_f] \rightarrow \Sph^2$ be a differentiable curve for $i\in\{1,\ldots,n\}$. The variation is a parameterized curve defined as $q_i^{\epsilon}: (-c,c)\times [t_0,t_f]\rightarrow\Sph^2$ for $c>0$, such that $q_i^0(t)= q_i(t)$ for any $t\in[t_0,t_f]$ and $q_i^{\epsilon}(t_0)=q_i(t_0)$, $q_i^{\epsilon}(t_f)=q(t_f)$ for any $\epsilon\in(-c,c)$.  

We can express the variation of the curve $q_i$ using the matrix exponential map as follows:
\begin{align}
    q_i^{\epsilon}(t) = e^ {\epsilon S(\gamma_i(t))} \,q_i(t),\label{eqn:qie}
\end{align}
for differentiable curves $\gamma_i: [t_0,t_f] \rightarrow \Re^3$ satisfying $\gamma_i(t_0)=\gamma_i(t_f)=0$.  Since the exponent $\epsilon S(\gamma_i)$ is skew-symmetric, the exponential matrix is in the special orthogonal group, $\SO=\{R\in\Re^{3\times 3}\,|\, R^T R=I,\,\mathrm{det}[R]=1\}$, thereby guaranteeing that the variation is a parameterized curve on $\Sph^2$, i.e., $\|q_i^{\epsilon}(t)\|=1$. There is no loss of generality in requiring that $\gamma_i(t) \cdot q_i(t) = 0$ for all $t_0 \le t \le t_f$; that is, $\gamma_i$ and $q_i$ are orthogonal. In short, the variation of the curve $q_i$ in $\Sph^2$ is expressed in terms of a curve in $\Re^3$ via \refeqn{qie}.

The corresponding infinitesimal variations\index{variation!infinitesimal} are given by
\begin{align}
    \delta q_i(t) & = \frac{d}{d\epsilon}\bigg|_{\epsilon=0} q^{\epsilon}(t)  = S(\gamma_i(t)) q_i(t),\label{eqn:delqi5}
\end{align}
and satisfy $\delta q_i(t_0)=\delta q_i(t_f)=0$. Since the variation and the differentiation commute, the expression for the  infinitesimal variations of the time derivatives are given by
\begin{align}
    \delta \dot q_i(t) &= \frac{d}{d\epsilon}\bigg|_{\epsilon=0} \dot q^{\epsilon}(t) 
     =S(\dot \gamma_i(t)) q_i(t)+S(\gamma_i(t)) \dot q_i(t).\label{eqn:delqidot}
\end{align}
These expressions define the infinitesimal variations for a vector function $q=(q_1,\dots,q_n):[t_0,t_f] \rightarrow (\Sph^2)^n$.    The infinitesimal variations are important ingredients  to derive the Euler--Lagrange equations on $(\Sph^2)^n$.  We subsequently suppress the time argument, thereby simplifying the notation. 



The action integral\index{action integral} is the integral of the Lagrangian function\index{Lagrangian function} along a motion of the system over a fixed time period, i.e., $\mathfrak{G}=\int_{t_0}^{t_f} L(q,\dot q)\,dt$. The infinitesimal variation of the action integral is given by
\begin{align*}
\delta \mathfrak{G} =\frac{d}{d\epsilon}\bigg|_{\epsilon=0} \int_{t_0}^{t_f} L(q^\epsilon,\dot q^\epsilon)\, dt =0.
\end{align*}
This can be rewritten in terms of $\delta q_i$ as
\begin{align}
    \delta \mathfrak{G} =  \int_{t_0}^{t_f} \braces{ \sum_{i=1}^n  \deriv{L(q,\dot q)}{\dot q_i} \cdot \delta \dot q_i + \sum_{i=1}^n  \deriv{L(q,\dot q)}{q_i} \cdot \delta q_i}\,dt.
    \label{eqn:delG0}
\end{align}
Let $f_i\in\Re^3$ be the generalized force acting on $q_i$. The corresponding virtual work is given by
\begin{align}
\delta\mathcal{W} = \sum_{i=1}^n f_i\cdot \delta q_i.
\end{align}
According to the Lagrange--d'Alembert principle, we have $\delta\mathfrak{G}=-\delta\mathcal{W}$ for any variations. We now substitute the expressions for the infinitesimal variations of the motion \refeqn{delqi5} and \refeqn{delqidot} into this, and we simplify the result to obtain the Euler--Lagrange equations expressed in terms of $(q,\dot q)$.   

\begin{prop}
Consider a Lagrangian $L(q,\dot q):\T(\Sph^2)^n\rightarrow \Re$ for a mechanical system evolving on the product of two-spheres, with the generalized force $f_i$. The Euler-Lagrange equations are given by
\begin{align}
    \parenth{I_3-q_i q_i^T} \braces{ \frac{d}{dt}\parenth{\deriv{L(q,\dot q)}{\dot q_i}} - \deriv{L(q,\dot q)}{q_i}-f_i} = 0,\label{eqn:EL5b}
 \end{align}
for $i\in\{1,\ldots,n\}$. Here, the $3\times 3$ identity matrix is denoted by $I_3\in\Re^{3\times 3}$.
\end{prop}
\begin{proof}
See Appendix \ref{sec:Pf1}.
\end{proof}
This describes the evolution of the dynamical flow $(q,\dot q) \in \T(\Sph^2)^n$ on the tangent bundle\index{tangent bundle} of the configuration manifold\index{manifold!configuration} $(\Sph^2)^n$.
 In the above equation \refeqn{EL5b}, the expression in the braces corresponds to the Euler--Lagrange equation for dynamical systems evolving on $\Re^n$. It is interesting to note that the Euler--Lagrange equations on $\Sph^2$ corresponds to its orthogonal projection onto the plane normal to $q_i$, represented by the matrix $I_{3}-q_iq_i^T$.

Next, we consider the important case that the kinetic energy is a quadratic function of $\dot q$, and the potential energy $U$ is given as a function of $q$, i.e., the Lagrangian can be written as
 \begin{align}
    L(q,\dot q)=\frac{1}{2} \sum_{j,k=1}^n \dot q_j^T m_{jk}(q)\dot q_k - U(q), \label{eqn:QL5}
\end{align}
where the scalar inertial terms $m_{jk}: (\Sph^2)^n \rightarrow \Re^1$ satisfy the symmetry condition $m_{jk}(q)=m_{kj}(q)$ and the quadratic form in the time derivatives of the configuration is positive-definite on $(\Sph^2)^n$. 

\begin{cor}
The Euler--Lagrange equations for the Lagrangian given by \refeqn{QL5} with the generalized force $f_i$ are 
\begin{multline}
m_{ii}(q) \ddot q_i +(I_3-q_iq_i^T) \sum_{\substack{j=1\\j\neq i}}^n m_{ij}(q) \ddot q_j +m_{ii}(q)\norm{\dot q_i}^2q_i \\ +(I_3-q_iq_i^T)\braces{F_i(q,\dot q)+\deriv{U(q)}{q_i}-f_i}=0, \label{eqn:EL5c}
\end{multline}
for $i\in\{1,\ldots,n\}$, where $F_i(q,\dot q)\in\Re^3$ is
\begin{align*}
F_i(q,\dot q)= \sum_{j=1}^n \dot m_{ij}(q) \dot q_j- \frac{1}{2} \deriv{}{q_i} \sum_{j,k=1}^n \dot q_j^T m_{jk}(q) \dot q_k.
\end{align*}
\end{cor}
\begin{proof}
See Appendix \ref{sec:Pf1}.
\end{proof}
Note that the third term of \refeqn{EL5c} represents the centripetal acceleration to ensure that the unit-length constraint $\|q_i\|=1$ is always satisfied. If the inertia terms $m_{ij}$ are independent of the configuration $q$, then $F_i(q,\dot q)=0$.

\subsection{Euler--Lagrange equations in terms of $(q,\omega)$}

An alternate expression for the Euler--Lagrange equations is now obtained in terms of the angular velocities introduced in \refeqn{w}. We express the action integral in terms of the modified Lagrangian function\index{Lagrangian function!modified} 
\begin{align}
\tilde L(q,\omega)=L(q,S(\omega) q),\label{eqn:Lqw}
\end{align} 
where the kinematics equations are given by \refeqn{rotkin5}. We use the notation $\omega = (\omega_1,\dots,\omega_n) \in \Re^{3n}$, and we view the modified Lagrangian $\tilde L(q,\omega)$ as being defined on the tangent bundle $\T(\Sph^2)^n$. 

Let the modified action integral be $\tilde{\mathfrak{G}}=\int_{t_0}^{t_f}\tilde L\,dt$. Its infinitesimal variation  can be written as
\begin{align}
    \delta \tilde{\mathfrak{G}} =  \int_{t_0}^{t_f}  \sum_{i=1}^n \braces{ \deriv{\tilde L(q,\omega)}{\omega_i} \cdot \delta \omega_i + \deriv{\tilde L(q,\omega)}{q_i} \cdot \delta q_i}\,dt.\label{eqn:delG1}
\end{align}

Next, we derive expressions for the infinitesimal variation of the angular velocity vectors.  From \refeqn{w},
\begin{align*}
\delta\omega_i = S(\delta q_i) \dot q_i + S(q_i) \delta\dot q_i,
\end{align*}
Substituting \refeqn{delqi5} and \refeqn{delqidot} and rearranging,
\begin{align*}
\delta\omega_i = (S(\gamma_i) q_i)\times \dot q_i + S(q_i) (S(\dot\gamma_i) q_i +S(\gamma_i) \dot q_i).
\end{align*}
Expanding each term and using the fact that $q_i\cdot q_i=1$ and $q_i\cdot \dot q_i =q_i\cdot\gamma_i=0$, this reduces to
\begin{align*}
\delta\omega_i &= (\dot q_i\cdot\gamma_i) q_i +\dot\gamma_i -(q_i\cdot\dot\gamma_i)q_i.
\end{align*}
Substitute \refeqn{rotkin5} to obtain
\begin{align*}
\delta\omega_i &= (\gamma_i\cdot(S(\omega_i) q_i)) q_i + (I_3-q_iq_i^T)\dot\gamma_i,\\
  &= q_iq_i^T (S(\gamma_i)\omega_i) + (I_3-q_iq_i^T)\dot\gamma_i.
\end{align*}
The matrix $q_iq_i^T$ corresponds to the orthogonal projection along $q_i$. But, as both of $\gamma_i$ and $\omega_i$ are orthogonal to $q_i$, $S(\gamma_i)\omega_i$ is already parallel to $q_i$. Therefore, 
\begin{align}
\delta\omega_i = -S(\omega_i) \gamma_i + (I_3-q_iq_i^T)\dot\gamma_i.\label{eqn:avv5}
\end{align}
The infinitesimal variation of $\omega_i$ is composed of two parts: the first term $-S(\omega_i) \gamma_i=\gamma_i\times\omega_i$ is parallel to $q_i$, and it represents the variations due to the change of $q_i$; the second term corresponds to the orthogonal projection of $\dot\gamma_i$ onto the orthogonal complement to $q_i$, and it is due to the time rate change of the variation of $q_i$. 

We now substitute \refeqn{delqi5} and \refeqn{avv5} into \refeqn{delG1}, and simplify the result to obtain the Euler--Lagrange equations expressed in terms of $(q,\omega)$.   

\begin{prop}
The Euler--Lagrange equations on $(\Sph^2)^n$ for the Lagrangian given by \refeqn{Lqw} with the generalized force $f_i$ are
\begin{multline}
 (I-q_iq_i^T) \braces{ \frac{d}{dt} \parenth{\deriv{\tilde L(q,\omega)}{\omega_i}} -2S(\omega_i) \deriv{\tilde L(q,\omega)}{\omega_i}} \\
 -   S(q_i)\parenth{ \deriv{\tilde L(q,\omega)}{q_i} +f_i}=0,\quad i=1,\dots,n. \label{eqn:ELom5}
\end{multline}
\end{prop}
\begin{proof}
See Appendix \ref{sec:Pf2}.
\end{proof}

This form of the Euler--Lagrange equations on $(\Sph^2)^n$, expressed in terms of angular velocities, can be obtained directly from the Euler--Lagrange equations given in \refeqn{EL5b} by viewing the kinematics \refeqn{rotkin5} as defining a change of variables from $\dot q$ to $\omega$.  This establishes the equivalence of the Euler--Lagrange equations in terms of $(q,\omega)$ \refeqn{ELom5} and the Euler--Lagrange equations in terms of $(q,\dot q)$ \refeqn{EL5b}.

Next, we consider the important case that the kinetic energy is a quadratic form as in \refeqn{QL5}. Substituting \refeqn{rotkin5} into \refeqn{QL5}, and using the fact that $\omega_i^T S(q_i)^TS(q_i) \omega_i = \omega_i^T \parenth{I_3 - q_iq_i^T} \omega_i = \omega_i^T \omega_i$, the modified Lagrangian can be expressed as
\begin{multline}
\tilde L(q,\omega)=\frac{1}{2} \sum_{i=1}^n  \omega_i^T m_{ii}(q) \omega_i \\
+ \frac{1}{2} \sum_{i=1}^n  \sum_{\substack{j=1\\j \ne i}}^n  \omega_i^T S(q_i)^T m_{ij}(q) S(q_j) \omega_j - U(q).\label{eqn:QLqw}
\end{multline}
Substituting this into \refeqn{ELom5} yields the corresponding Euler--Lagrange equations as follows.

\begin{cor}
The Euler--Lagrange equations for the modified Lagrangian given by \refeqn{QLqw} with the generalized force $f_i$ are
\begin{multline}
m_{ii}(q) \dot \omega_i +  \sum_{\substack{j=1\\j\neq i}}^n S(q_i)^T m_{ij}(q) S(q_j)\dot \omega_j \\
-  m_{ij}(q)S(q_i)\norm{\omega_j}^2 q_j \\
+ S(q_i) \braces{F_i(q,\omega) +  \deriv{U(q)}{q_i}-f_i}=0, \label{eqn:ELav55}
\end{multline}
for $i=1,\ldots,n$, where $F_i(q,\omega)\in\Re^3$ is
\begin{align*}
F_i(q,\omega)&= \sum_{j=1}^n \dot m_{ij}(q) S(\omega_j) q_j\\ &\quad - \frac{1}{2} \sum_{j=1}^n  \sum_{k=1}^n   (q_j^T S(\omega_j)^T  S(\omega_k) q_k) \deriv{m_{jk}(q)}{q_i}.
\end{align*}
\end{cor}
\begin{proof}
See Appendix \ref{sec:Pf2}.
\end{proof}
Similar to Corollary 1, if the inertial terms are independent of the configuration, then $F_i(q,\omega)=0$. This version of the Euler--Lagrange differential equations describe the dynamical flow $(q,\omega) \in \T(\Sph^2)^n$ on the tangent bundle\index{tangent bundle} of $(\Sph^2)^n$.

\section{Hamiltonian Mechanics on Two-Spheres}

We will now introduce the Legendre transformation and then we derive Hamilton's equations for dynamics that evolve on $(\Sph^2)^n$.   The derivation is based on the phase space variational principle, a natural modification of Hamilton's principle for Lagrangian dynamics.   Two forms of Hamilton's equations are obtained.  One form is expressed in terms of momentum vectors $(\mu_1,\dots,\mu_n) \in \T^*_q (\Sph^2)^n$ that are conjugate to the velocities $(\dot q_1,\dots,\dot q_n) \in \T_q (\Sph^2)^n$, where $q \in (\Sph^2)^n$.   The other form of Hamilton's equations are expressed in terms of momentum vectors $(\pi_1,\dots,\pi_n) \in \Re^{3n}$ that are conjugate to the angular velocities $(\omega_1,\dots,\omega_n) \in \Re^{3n}$.

\subsection{Hamilton's equations in terms of $(q,\mu)$}\label{subsec:LTTS}

As in the prior section, we begin with a Lagrangian function\index{Lagrangian function} $L:\T(\Sph^2)^n\rightarrow\Re^1$, which is a real-valued function defined on the tangent bundle of the configuration manifold\index{manifold!configuration} $(\Sph^2)^n$. The Legendre transformation\index{Legendre transformation} of the Lagrangian function\index{Lagrangian function} $L(q,\dot q)$ leads to the Hamiltonian form of the equations of motion in terms of a conjugate momentum\index{conjugate momentum} vector. For $q_i\in\Sph^2$, the corresponding conjugate momentum $\mu_i$ lies in the dual space $\T_{q_i}^* \Sph^2$. We identify the tangent space $\T_{q_i} \Sph^2$ and its dual space $\T_{q_i}^* \Sph^2$ by using the usual dot product in $\Re^3$. More explicitly, the Legendre transformation is given by
\begin{align*}
    \mu_i \cdot \dot q_i & = \deriv{ L(q,\dot q)}{\dot q_i}\cdot \dot q_i,
\end{align*}
for any $\dot q_i \in \Re^3$ orthogonal to $q_i$. Since the component of $\mu_i$ parallel to $q_i$ has no effect on the inner product above, the vector representing $\mu_i$ is selected to be orthogonal to $q_i$; that is $\mu_i$ is equal to the projection of $\deriv{L(q,\dot q)}{\dot q_i}$ onto the tangent space $\T_{q_i} \Sph^2$.  Thus
\begin{align}
    \mu_i =  (I_3-q_i q_i^T)  \deriv{L(q,\dot q)}{\dot q_i}.\label{eqn:P2SHaa}
\end{align}


We assume that the Lagrangian function\index{Lagrangian function} has the property that the Legendre transformation is invertible in the sense that the above $n$ algebraic equations, viewed as a mapping from $\T_q(\Sph^2)^n$ to $\T^*_q(\Sph^2)^n$, is invertible.   Since these tangent and cotangent spaces are embedded in $\Re^{3n}$, we can view the Legendre transformation as being the restriction of a mapping from $\Re^{3n}$ to $\Re^{3n}$ that is invertible. 

The Hamiltonian function\index{Hamiltonian function} $H: \T^*(\Sph^2)^n \rightarrow \Re^1$ is given by 
\begin{align}
H(q,\mu)=\sum_{i=1}^{n} \mu_i \cdot \dot q_i-L(q,\dot q), \label{eqn:5Hom}
\end{align}
where the right hand side is expressed in terms of $(q,\mu)$ using the Legendre transformation \refeqn{P2SHaa}.  

The Legendre transformation can be viewed as defining a transformation $(q,\dot q) \mapsto (q,\mu)$, which implies that the Euler--Lagrange equations can be written in terms of the transformed variables; this is effectively Hamilton's equations.   However, Hamilton's equations can also be obtained using Hamilton's phase space variational principle, and this approach is now introduced.  

Consider the action integral in the form,
\begin{align*}
\mathfrak{G}=\int_{t_0}^{t_f}\braces{ \sum_{i=1}^n\mu_i\cdot \dot q_i - H(q,\mu)}\, dt.
\end{align*}
Integrating by parts and using the fact that the variation $\delta q$ vanishes at $t_0$ and $t_f$, the infinitesimal variation of the action integral is given by
\begin{align}
\delta \mathfrak{G} & =
\sum_{i=1}^n\int_{t_0}^{t_f}\Bigg\{\parenth{-\dot \mu_i -\deriv{H(q,\mu)}{q_i}}\cdot \delta q_i\nonumber\\&\quad+\parenth{\dot q_i-\deriv{H(\mu,p)}{\mu_i}} \cdot \delta \mu_i\Bigg\}\,dt = 0.
\label{eqn:delG2}
\end{align}

Next, we derive the expression for the variations of $\mu_i$. According to the definition of the conjugate momenta $\mu_i$ given by \refeqn{P2SHaa}, we have $q_i\cdot \mu_i=0$, which implies that $\delta q_i \cdot \mu_i + q_i\cdot\delta\mu_i=0$.   To impose this constraint on the variations explicitly, we decompose $\delta\mu_i$ into the sum of two orthogonal components: one component parallel to $q_i$, namely $\delta\mu_i^C=q_iq_i^T\delta\mu_i$, and the other component orthogonal to $q_i$, namely $\delta\mu_i^M=(I_{3\times 3}-q_iq_i^T)\delta\mu_i$.   Satisfaction of the constraint implies that $q_i^T\delta\mu_i^C = q_i^T\delta\mu_i= -\mu_i^T\delta q_i$, so that $\delta \mu_i^M = (I_3-q_iq_i^T) \delta \mu_i$ is otherwise unconstrained. Substituting this and \refeqn{delqi5} into \refeqn{delG2}, we obtain Hamilton's equations in terms of $(q,\mu)$ as follows.

\begin{prop}
Hamilton's equations on $(\Sph^2)^n$ for the Hamiltonian given by \refeqn{5Hom} with the generalized force $f_i$ are 
\begin{align}
\dot q_i&=(I_{3\times 3}-q_i q_i^T)\deriv{H(q,\mu)}{\mu_i}, \label{eqn:5HEmua}\\
\dot \mu_i&= -(I_{3\times 3}-q_i q_i^T)\parenth{\deriv{H(q,\mu)}{q_i}-f_i}\nonumber\\
&\quad +\deriv{H(q,\mu)}{\mu_i}\times \parenth{\mu_i\times q_i},\label{eqn:5HEmub}
\end{align}
for $i=1,\ldots,n$.
\end{prop}
\begin{proof}
See Appendix \ref{sec:Pf3}.
\end{proof}
Thus, equations \refeqn{5HEmua} and \refeqn{5HEmub} describe the Hamiltonian flow in terms of $(q,\mu) \in \T^*(\Sph^2)^n$ on the cotangent bundle\index{cotangent bundle} of $(\Sph^2)^n$.

When $f_i=0$, any time-independent Hamiltonian is preserved along the solution of Hamilton's equations, since
\begin{align*}
\frac{dH}{dt}& = \deriv{H}{t} +\sum_{i=1}^n\deriv{H}{q_i}\cdot \dot q_i + \deriv{H}{\mu_i}\cdot \dot \mu_i\\
& = \deriv{H}{t} +\sum_{i=1}^n 
\deriv{H}{\mu_i} \cdot \braces{\deriv{H}{\mu_i}\times (\mu_i\times q_i)} = \deriv{H}{t}.
\end{align*}

Next, we consider the case where the kinetic energy is a quadratic function of the time derivatives of the configuration so that the Lagrangian is given by \refeqn{QL5}. The conjugate momentum\index{conjugate momentum} vector is defined by the Legendre transformation
\begin{align*}
    \mu_i  
    &=m_{ii}(q) \dot q_i + (I_3-q_i q_i^T) \sum_{\substack{j=1\\j\neq i}}^n  m_{ij}(q)\dot q_j.
\end{align*}
We assume that these algebraic equations, viewed as a linear mapping from $(\dot q_1,\dots,\dot q_n) \in \T_q (\Sph^2)^n \subset \Re^{3n}$ to $(\mu_1,\dots,\mu_n) \in \T^*_q (\Sph^2)^n \subset \Re^{3n}$, can be inverted and expressed in the form
\begin{align}
\dot q_i = (I_3-q_iq_i^T) \sum_{j=1}^n m^I_{ij}(q) \mu_j,\label{eqn:5HEpa}
\end{align}
where $m^I_{ij}:(\Sph^2)^n \rightarrow \Re^{3 \times 3}$. There is no loss of generality in including the indicated projection in the above expression since the inverse necessarily guarantees that if $(\mu_1,\dots,\mu_n) \in \T^*_q (\Sph^2)^n \subset \Re^{3n}$ then 
$(\dot q_1,\dots,\dot q_n) \in \T_q (\Sph^2)^n \subset \Re^{3n}$. The Hamiltonian \index{Hamiltonian function} can be expressed as
\begin{align}
H(q,\mu)= \frac{1}{2} \sum_{j,k=1}^n \mu_j^T m^I_{jk}(q) \mu_k + U(q).\label{eqn:QHqmu}
\end{align}

\begin{cor}
Hamilton's equations for the Hamiltonian given by \refeqn{QHqmu} with the generalized force $f_i$ are \refeqn{5HEpa} and 
\begin{align}
    \dot \mu_i = & \sum_{j=1}^n \parenth{m^I_{ij}(q)\mu_j} \times \parenth{\mu_i \times q_i}\nonumber\\
    &  - (I_3-q_i q_i^T) \frac{1}{2} \deriv{}{q_i}  \sum_{j=1}^n \sum_{k=1}^n \mu_j^T m^I_{jk}(q) \mu_k \nonumber \\
    &-(I_3-q_i q_i^T) \parenth{\deriv{U(q)}{q_i}-f_i}.   \label{eqn:5HEpb}
\end{align}
\end{cor}
\begin{proof}
See Appendix \ref{sec:Pf3}.
\end{proof}

Hamilton's equations \refeqn{5HEpa} and \refeqn{5HEpb}
describe the Hamiltonian flow in terms of $(q,\mu) \in \T^*(\Sph^2)^n$ on the cotangent bundle\index{cotangent bundle} of $(\Sph^2)^n$. 


\subsection{Hamilton's equations in terms of $(q,\pi)$}\label{subsec:LTTS5b}

We now present an alternate version of Hamilton's equations using the Legendre transformation of the modified Lagrangian function\index{Lagrangian function} $\tilde L(q,\omega)$ to define the conjugate momentum\index{conjugate momentum} vector.   The Legendre transformation\index{Legendre transformation} $(\omega_1,\dots,\omega_n) \in \Re^{3n} \rightarrow (\pi_1,\dots,\pi_n) \in \Re^{3n}$ is defined by
\begin{align}
    \pi_i = (I_3 -q_iq_i^T) \deriv{\tilde L(q,\omega)}{\omega_i}.  \label{eqn:P2SHa}
\end{align}
Here $\pi_i \in \Re^3$ is viewed as conjugate to $\omega_i \in \Re^3$, $i=1,\dots,n$.   We use the notation $\pi=(\pi_1,\dots,\pi_n) \in \Re^{3n}$.   
We assume that the modified Lagrangian function\index{Lagrangian function} has the property that the Legendre transformation is invertible in the sense that the above algebraic equations, viewed as a mapping from $\Re^{3n}$ to $\Re^{3n}$, is invertible.   

The modified Hamiltonian function\index{Hamiltonian function!modified} is given by 
\begin{align}
\tilde H(q,\pi)=\sum_{j=1}^{n} \pi_j \cdot \omega_j-\tilde L(q,\omega),
\label{eqn:5Hompi}
\end{align}
where the right hand side is expressed in terms of $(q,\pi)$ using the Legendre transformation \refeqn{P2SHa}.

Consider the modified action integral of the form,
\begin{align*}
\tilde{\mathfrak{G}} = \int_{t_0}^{t_f} \Bigg\{ \sum_{j=1}^n\pi_j \cdot \omega_j - \tilde{H}(q,\omega)\Bigg\}\,dt.
\end{align*}
Take the infinitesimal variation of $\tilde{\mathfrak{G}}$ and integrate by parts to obtain
\begin{align}
\delta\tilde{\mathfrak{G}} =
& \sum_{j=1}^n \int_{t_0}^{t_f}  \parenth{ \omega_i-\deriv{\tilde H(q,\pi)}{\pi_i} }\cdot\delta\pi_i\nonumber\\ &+ \parenth{ -\dot\pi_i+S(\omega)\pi_i -S(q_i)\deriv{\tilde H(q,\pi)}{q_i}}\cdot\gamma_i\, dt,\label{eqn:delG3}
\end{align}
where we use the fact that $(I_{3\times 3}-q_iq_i^T)\pi_i = \pi_i$ since $\pi_i$ is orthogonal to $q_i$ by the definition~\refeqn{P2SHa}. 

The orthogonality condition $\pi_i \cdot q_i=0$ also implies that $\delta q_i\cdot\pi_i +q_i\cdot\delta\pi_i=0$.  To impose this constraint on the variations explicitly, we decompose $\delta\pi_i$ into a component  that is parallel to $q_i$, namely $\delta\pi_i^C = q_iq_i^T\delta\pi_i$, and a component that is orthogonal to $q_i$, namely $\delta\pi_i^M =(I_{3\times 3}-q_iq_i^T)\delta\pi_i$.      From the above constraint, we have $q_i^T\delta\pi_i = -\pi_i^T\delta q_i = -\pi_i^T S(\gamma_i)q_i = \pi_i^T S(q_i)\gamma_i$.    Therefore $\delta\pi_i^C = q_iq_i^T\delta\pi_i = q_i\pi_i^T S(q_i)\gamma_i$.

\begin{prop}
Hamilton's equations for the modified Hamiltonian given by \refeqn{5Hompi} with the generalized force $f_i$ are
\begin{align}
\dot q_i & = -S(q_i)\deriv{\tilde H(q,\pi)}{\pi_i}, \label{eqn:5HEpifa} \\
\dot\pi_i & = -S(q_i)\deriv{\tilde H(q,\pi)}{q_i}+\deriv{\tilde H(q,\pi)}{\pi_i}\times \pi_i+S(q_i)f_i, \label{eqn:5HEpifb}
\end{align}
for $i=1,\ldots,n$.
\end{prop}
\begin{proof}
See Appendix \ref{sec:Pf4}.
\end{proof}
Thus equations \refeqn{5HEpifa} and \refeqn{5HEpifb} describe the Hamiltonian flow in terms of $(q,\pi) \in \T^*(\Sph^2)^n$ on the cotangent bundle\index{cotangent bundle} of $(\Sph^2)^n$. 

When $f_i=0$, any time-independent modified Hamiltonian is preserved along the solution of Hamilton's equations, since
\begin{align*}
\frac{d\tilde H}{dt}& = \deriv{\tilde H}{t} +\sum_{i=1}^n\deriv{\tilde H}{q_i}\cdot \dot q_i + \deriv{\tilde H}{\pi_i}\cdot \dot \pi_i\\
& = \deriv{\tilde H}{t} +\sum_{i=1}^n 
\deriv{\tilde H}{\pi_i} \cdot \braces{\deriv{\tilde H}{\pi_i}\times \pi} = \deriv{\tilde H}{t}.
\end{align*}

We now consider the important case where the kinetic energy is a quadratic\index{kinetic energy!quadratic} function of the angular velocities in the form that arises from the Lagrangian given by \refeqn{Lqw}. The conjugate momentum is defined by the Legendre transformation
\begin{align}
    \pi_i  
    = & m_{ii}(q) \omega_i + \sum_{\substack{j=1\\j\neq i}}^n  S(q_i)^T m_{ij}(q) S(q_j) \omega_j. \label{eqn:HEmomm5}
\end{align}
We assume these algebraic equations, viewed as a linear mapping from $(\omega_1,\dots,\omega_n) \in \Re^{3n}$ to $(\pi_1,\dots,\pi_n) \in \Re^{3n}$ can be inverted and expressed in the form
\begin{align}
\omega_i = \sum_{j=1}^n m^I_{ij}(q) \pi_j,\label{eqn:wi3}
\end{align}
where $m^I_{ij}:(\Sph^2)^n \rightarrow \Re^{3 \times 3}$. The modified Hamiltonian function\index{Hamiltonian function!modified} can be expressed as
\begin{align}
\tilde H(q,\pi)= \frac{1}{2} \sum_{i=1}^n \sum_{j=1}^n  \pi_i^T m^I_{ij}(q) \pi_j + U(q). \label{eqn:5Hpi}
\end{align}

\begin{cor}
Hamilton's equations for the modified Hamiltonian given by \refeqn{5Hpi} are 
\begin{align}
\dot q_i & = -S(q_i) \braces{\sum_{j=1}^n m^I_{ij}(q) \pi_j}, \label{eqn:5HEpigg} \\
\dot \pi_i & = -S(q_i) \braces{\frac{1}{2} \deriv{}{q_i} \sum_{j,k=1}^n \pi_j^T m^I_{jk}(q)  \pi_k + \deriv{U(q)}{q_i}}\nonumber \\
&\quad + \braces{\sum_{j=1}^n m^I_{ij}(q) \pi_j} \times \pi_i+S(q_i)f_i, \label{eqn:5HEpihh}
\end{align}
for $i=1,\ldots,n$.
\end{cor}
\begin{proof}
See Appendix \ref{sec:Pf4}.
\end{proof}

\section{Dynamics on Chain Pendulum}

A chain pendulum is a connection of $n$ rigid links, that are serially connected by two degree-of-freedom spherical joints.  We assume that each link of the chain pendulum is a rigid rod with mass concentrated at the outboard end of the link.   One end of the chain pendulum is connected to a spherical joint that is supported by a fixed base.   A constant gravitational acceleration acts on each link of the chain pendulum.  This may represent a spherical pendulum ($n=1$), or a double spherical pendulum ($n=2$) as special cases. 

We demonstrate that globally valid Euler--Lagrange equations can be developed for the chain pendulum, and they can be expressed in a compact form.   The results provide an intrinsic and unified framework to study the dynamics of a chain pendulum system, that is applicable for an arbitrary number of links, and globally valid for any configuration of the links.

The mass of the $i$-th link is denoted by $m_i$ and the link length is denoted by $l_i$.  For simplicity, we assume that the mass of each link is concentrated at the outboard end of the link. An inertial frame is chosen such that the first two axes are horizontal and the third axis is vertical.   The origin of the inertial frame is located at the fixed spherical joint.    Each of the chain links has a body-fixed frame with the third axis of the body-fixed frame aligned with the axial direction of the link.     The vector $q_1\in\Sph^2$ represents the direction from the fixed base to the mass element of the first link, and the vector $q_i\in\Sph^2$ represents the direction from the $(i-1)$-th spherical joint to the concentrated mass element of the $i$-th link.  Thus, the configuration of the chain pendulum is the ordered $n$-tuple of configurations of each link $q=(q_1,\dots,q_n) \in (\Sph^2)^n$, so that the configuration manifold\index{manifold!configuration} is $(\Sph^2)^n$.   The chain pendulum has $2n$ degrees of freedom. 

Let $x_i\in\Re^3$ be the position of the outboard end of the $i$-th link in the inertial frame; it can be written as $x_i =  \sum_{j=1}^i l_jq_j$. The total kinetic energy is composed of the kinetic energy of each mass:
\begin{align*}
T(q,\dot q) 
 =\frac{1}{2}\sum_{i=1}^n m_i \| \sum_{j=1}^i l_j\dot q_j\|^2.
 \end{align*}
 This can be rewritten as
\begin{align}
T(q,\omega) = \frac{1}{2} \sum_{i,j=1}^n M_{ij}l_il_j  \omega_i^TS(q_i)^TM_{ij}l_il_j S(q_j)\omega_j,\label{eqn:Tcp}
\end{align}
where the real inertia constants $M_{ij}$ are given by
\begin{align*}
M_{ij} =  \parenth{\sum_{k=\max\{i,j\}}^n m_k},\quad i,j=1,\dots,n. 
\end{align*}

The potential energy consists of the gravitational potential energy of all mass elements.  The potential energy can be written as
\begin{align}
U(q) = \sum_{i=1}^n m_i g e_3^T x_i  =  \sum_{i=1}^n \sum_{j=i}^n m_jg l_i e_3^T q_i.
\label{eqn:Ucp}
\end{align}
The modified Lagrangian function\index{Lagrangian function} $\tilde L:\T (\Sph^2)^n \rightarrow \Re^1$ of the chain pendulum is given by $\tilde L(q,\omega)=T(q,\omega)-U(q)$ from \refeqn{Tcp} and \refeqn{Ucp}. 

Also, suppose that there exist a control torque $\tau\in\Re^3$ acting on the spherical joint connecting the fixed base and the first link, and a disturbance force $d\in\Re^3$ acting on the tip of the last link. The corresponding virtual work is given by
\begin{align*}
\delta\mathcal{W} 
=\gamma_1 \cdot \tau + \sum_{i=1}^n l_i S(\gamma_i)q_i \cdot d
=\gamma_1 \cdot \tau + \sum_{i=1}^n \gamma_i \cdot l_iS(q_i)d.
\end{align*}
Therefore, the generalized forces are given by $f_1=\tau + l_1S(q_1)d$, and $f_j=l_jS(q_j)d$ for $j\geq 2$.

Substituting this into \refeqn{ELav55}, the Euler--Lagrange equations for a chain pendulum are given by
\begin{align}
& M_{ii} l_i^2 \dot \omega_i + \sum_{\substack{j=1\\j \ne i}}^n M_{ij}l_il_j S^T(q_i)S(q_j) \dot \omega_j\nonumber \\ & - \sum_{\substack{j=1\\j \ne i}}^n M_{ij}l_il_j \norm{\omega_j}^2 S(q_i)q_j - \sum_{j=i}^n m_j g l_i S(e_3) q_i=S(q_i)f_i.
\end{align}

Similarly, the Legendre transformation is given by \refeqn{HEmomm5}, and from \refeqn{5HEpigg}, \refeqn{5HEpihh}, Hamilton's equations can be written as
\begin{align}
\dot q_i & = -S(q_i) \braces{\sum_{j=1}^n M^I_{ij}(q) \pi_j}, \\
\dot \pi_i & = -S(q_i) \braces{\frac{1}{2} \deriv{}{q_i} \sum_{j,k=1}^n \pi_j^T M^I_{jk}(q)  \pi_k + \sum_{j=1}^n m_j gl_j e_3}\nonumber \\
&\quad + \braces{\sum_{j=1}^n M^I_{ij}(q) \pi_j} \times \pi_i+S(q_i)f_i.
\end{align}

These are remarkably compact considering the complexity of the dynamics, and they are well structured compared with the equations of motion expressed in terms of angles. 

This mathematical model may be applied to a wide class of other dynamical systems, such as articulated robotic systems. As they are developed for an arbitrary number of links, they are readily extended to finite-element approximations of cables or slender rods after augmenting the potential with an elastic potential term. The proposed global formulations avoid singularities associated with local coordinates.

\appendix

\subsection{Hat map}
Several properties of the hat map are summarized as follows.
\begin{gather}
    S(x) y = x\times y = - y\times x = - S(y) x,\\
    S(x)^2 = -(x^T x) I_{3\times 3} + x x^T,\label{eqn:S2}\\
    S(x)^3 = -(x^T x) S(x),\label{eqn:S3}\\
    x\cdot S(y)z = y\cdot S(z)x = z\cdot S(x)y,\label{eqn:STP}\\
    S(x) S(y) z=(x\cdot z) y - (x\cdot y)z=(yx^T-x^TyI_{3\times 3})z,\label{eqn:VTP}\\
    S(x\times y) = S(x) S(y) -S(y)S(x) = yx^T-xy^T,\label{eqn:Sxy}
\end{gather}
for any $x,y,z\in\Re^3$.

\subsection{Proof of Proposition 1}\label{sec:Pf1}
Substituting  \refeqn{delqi5} and \refeqn{delqidot} into \refeqn{delG0}, and rearranging with \refeqn{STP},
\begin{align*}
    \delta \mathfrak{G} & = \int_{t_0}^{t_f}  \sum_{i=1}^n \Bigg\{  \dot \gamma_i \cdot \parenth{S(q_i) \deriv{L(q,\dot q)}{\dot q_i}} \\
    &\quad + \gamma_i \cdot \parenth{S(\dot q_i) \deriv{L(q,\dot q)}{\dot q_i}+ S(q_i) \deriv{L(q,\dot q)}{q_i}}\Bigg\}\,dt.
\end{align*}
Integrating the first term on the right by parts, and using the fact that the variation vanishes at $t_0$ and $t_f$, it reduces to
\begin{multline*}
    \delta \mathfrak{G}  = - \sum_{i=1}^n \int_{t_0}^{t_f} \gamma_i \cdot \parenth{S(q_i) \braces{ \frac{d}{dt} \deriv{L(q,\dot q)}{\dot q_i} -  \deriv{L(q,\dot q)}{q_i}}}\, dt.
\end{multline*}
The virtual work is written as
\begin{align*}
\delta \mathcal{W} = \sum_{i=1}^n f_i\cdot S(\gamma_i)q_i = \sum_{i=1}^n \gamma_i\cdot S(q_i)f_i.
\end{align*}
According to Lagrange--d'Alembert principle, $\delta \mathfrak{G}=-\delta\mathcal{W}$ for all continuous variations $\gamma_i:[t_0,t_f] \rightarrow \Re^3,$ that satisfy $(\gamma_i \cdot q_i)=0$. The fundamental lemma of the calculus of variations implies that the expression $\frac{d}{dt}\deriv{L}{\dot q_i}-\deriv{L}{q_i}-f_i$ is parallel to $q_i$, or equivalently, \refeqn{EL5b}.

Next, consider the case where the kinetic energy is given as a quadratic form as \refeqn{QL5}. Substituting \refeqn{QL5} into \refeqn{EL5b},
\begin{multline}
(I_3-q_iq_i^T) \bigg\{ \sum_{j=1}^n m_{ij}(q)\ddot q_j +  \sum_{j=1}^n \dot m_{ij}(q) \dot q_j\\ - \frac{1}{2} \deriv{}{q_i} \sum_{j=1}^n \sum_{k=1}^n \dot q_j^T m_{jk}(q) \dot q_k +  \deriv{U(q)}{q_i}-f_i\bigg\}=0. \label{eqn:ELQ5b}
\end{multline}
Since $q_i^T \dot q_i = 0$, it follows that $\frac{d}{dt}(q_i^T \dot q_i)= (q_i^T \ddot q_i) + \norm{ \dot q_i}^2=0$; thus
we obtain
\begin{align*}
   (I_3 -q_iq_i^T)  \ddot q_i & = \ddot q_i - (q_i q_i^T) \ddot q_i =  \ddot q_i + \norm{\dot q_i}^2  q_i.
\end{align*}
Substituting this into \refeqn{ELQ5b} yields \refeqn{EL5c}.

\subsection{Proof of Proposition 2}\label{sec:Pf2}

Substituting \refeqn{delqi5} and \refeqn{avv5} into \refeqn{delG1} and using the integration by parts, 
\begin{align*}
    \delta  \tilde{ \mathfrak{G}} & = \sum_{i=1}^n \int_{t_0}^{t_f} \gamma_i \cdot \Bigg\{ -\frac{d}{dt} \parenth{(I-q_iq_i^T) \deriv{\tilde L(q,\omega)}{\omega_i}} \\
&\quad    +   S(q_i) \deriv{\tilde L(q,\omega)}{q_i}+S(\omega_i) \deriv{\tilde L(q,\omega)}{\omega_i} \Bigg\}  \, dt.
\end{align*}
According to Lagrange--d'Alembert principle, $\delta \tilde{\mathfrak{G}}=-\delta\mathcal{W}$ for all differentiable functions $\gamma_i: [t_0,t_f] \rightarrow  \Re^3$ that satisfy $(\gamma_i \cdot q_i)=0$ and vanish at $t_0$, and $t_f$. This implies that the expression in the braces of the equation below is parallel to $q_i$, or equivalently
\begin{multline}
S(q_i)^2 \Bigg\{ \frac{d}{dt} \parenth{S^2(q_i) \deriv{\tilde L(q,\omega)}{\omega_i}}
    +S(\omega_i) \deriv{\tilde L(q,\omega)}{\omega_i}\\
    +   S(q_i) \parenth{\deriv{\tilde L(q,\omega)}{q_i}+f_i}\Bigg\} =0,\label{eqn:ELom5_0}
\end{multline}
where we used $-S(q_i)^2=I_{3\times 3}-q_iq_i^T$ obtained by \refeqn{S2}.

To further simplify these expressions, we derive a few identities. From \refeqn{VTP} and $w_i\cdot q_i=0$, 
\begin{align*}
S(q_i)S(\omega_i)S(q_i) =S(q_i)\{-\omega_i^TqI_3 + q_i \omega_i^T \}=0.
\end{align*}
From this, \refeqn{rotkin5}, and \refeqn{Sxy}, it follows that
\begin{align*}
S(q_i)S(\dot q_i) &= S(q_i)\{ S(\omega_i) S(q_i) - S(q_i) S(\omega_i)\}\\
& = -S(q_i)^2S(\omega_i),
\end{align*}
and similarly, $S(\dot q_i) S(q_i)=S(\omega_i)S(q_i)^2$. Consequently, these results can be used to obtain
\begin{align*}
S(q_i)^2 \{ S(\dot q_i)S(q_i) +S(q_i)S(\dot q_i)\} = -S(q_i)^4 S(\omega_i).
\end{align*}
Substituting these into \refeqn{ELom5_0} and using \refeqn{S2} and \refeqn{S3} repeatedly yield \refeqn{ELom5}. 

Next, the Euler--Lagrange equations for the case that the kinetic energy is given as a quadratic form, namely \refeqn{ELav55} can be obtained by either substituting \refeqn{QLqw} into \refeqn{ELom5}, or rewriting \refeqn{EL5c} in terms of the angular velocity as a change of variables. Here, we follow the latter approach as the corresponding proof is more concise.

From the kinematics equation \refeqn{rotkin5} and \refeqn{VTP}, we have
\begin{align*}
\ddot q_i = \dot\omega_i \times q_i + \omega_i \times (\omega_i\times q_i) =
-S(q_i)\dot\omega_i -\|\omega_i\|^2 q_i.
\end{align*}
Substituting this into \refeqn{EL5c}, and rearranging it with $I_{3\times 3}-q_iq_i^T=-S(q_i)^2$, 
\begin{align*}
-S(q_i)\Big\{ m_{ii}(q)\dot\omega_i&+S(q_i)\sum_{\substack{j=1\\j\neq i}}^n 
m_{ij}(q)(-S(q_j)\dot\omega_j-\|\omega_j\|^2 q_j)\\
&+S(q_i) \Big(F_i(q,\omega) +\deriv{U(q)}{q_i}-f_i\Big)\Big\}=0.
\end{align*}
In the above equation, the left hand side becomes zero when the expressions in the braces are either zero or parallel to $q_i$. However, the second case is not possible as they are perpendicular to $q_i$ by the definition. This yields \refeqn{ELav55}.

\subsection{Proof of Proposition 3}\label{sec:Pf3}

Substituting $\delta\mu_i=\delta\mu_i^C+\delta\mu_i^M$ and \refeqn{delqi5} into \refeqn{delG2}, and rearranging it with $q_i^T\delta\mu_i^C=-\mu_i^T\delta q_i$,
\begin{align*}
& \delta \mathfrak{G}  =
\sum_{i=1}^n\int_{t_0}^{t_f}\Bigg\{ \parenth{-\dot \mu_i -\deriv{H(q,\mu)}{q_i}}\cdot \delta q_i\\
&\quad\quad +\parenth{q_iq_i^T\parenth{\dot q_i-\deriv{H(q,\mu)}{\mu_i}}} \cdot \delta \mu_i^C\\
&\quad\quad +\parenth{(I_{3\times 3}-q_i q_i^T)\parenth{\dot q_i-\deriv{H(q,\mu)}{\mu_i}}} \cdot \delta \mu_i^M\Bigg\}\,dt,\\
& = \sum_{i=1}^n\int_{t_0}^{t_f}\Bigg\{ S(q_i) \parenth{-\dot \mu_i -\deriv{H(q,\mu)}{q_i}
+\mu_i q_i^T\deriv{H(q,\mu)}{\mu_i}} \cdot  \gamma_i\\
&\quad\quad +\parenth{\dot q_i-(I_{3\times 3}-q_i q_i^T)\deriv{H(q,\mu)}{\mu_i}} \cdot \delta \mu_i^M\Bigg\}\,dt.
\end{align*}
We now invoke Hamilton's phase space variational principle that $\delta \mathfrak{G}=-\delta\mathcal{W}$ for all possible functions $\gamma_i:[t_0,t_f] \rightarrow \Re^3$ satisfying $\gamma_i \cdot q_i = 0$ and $\delta\mu_i^M:[t_0,t_f] \rightarrow \Re^3$ that are always orthogonal to $q_i$ for $i=1,\dots,n$. According to the fundamental lemma of the calculus of variations, 
\begin{gather*}
S(q_i) \parenth{\dot \mu_i +\deriv{H(q,\mu)}{q_i}
-\mu_i q_i^T\deriv{H(q,\mu)}{\mu_i}-f_i}=0.
\end{gather*}
We multiply this by $S(q_i)$ and use a matrix identity to obtain
\begin{gather*}
(I_{3\times 3}-q_i q_i^T)\parenth{\dot \mu_i +\deriv{H(q,\mu)}{q_i}
-\mu_i q_i^T\deriv{H(q,\mu)}{\mu_i}-f_i}=0.
\end{gather*}
Since both terms multiplying $\delta \mu_i^M$ in the above variational expression are necessarily orthogonal to $q_i$, it yields \refeqn{5HEmua}.

We now determine an expression for $\dot \mu_i$.    The above equation only determines the component of $\dot\mu_i$ that is normal to $q_i$. The other component of $\dot\mu_i$ that is parallel to $q_i$ is derived by taking the time derivative of $q_i\cdot \mu_i=0$ to obtain $q_i\cdot \dot\mu_i = -\dot q_i \cdot\mu_i$. Thus, $\dot\mu_i$ is obtained by the sum of two components as
\begin{align*}
\dot\mu_i & = (I_{3\times 3}-q_i q_i^T)\parenth{-\deriv{H}{q_i}
+\mu_i q_i^T\deriv{H}{\mu_i}} -(\mu_i^T\dot q_i) q_i,
\end{align*}
which is reduced to \refeqn{5HEmub} via \refeqn{5HEmua} and \refeqn{VTP}.

Next, \refeqn{5HEpb} can be derived by substituting \refeqn{QHqmu} into \refeqn{5HEmub} directly and rearrange it.

\subsection{Proof of Proposition 4}\label{sec:Pf4}

Using $\delta\pi_i=\delta\pi_i^C+\delta\pi_i^M$, $\delta\pi_i^C = q_iq_i^T\delta\pi_i = q_i\pi_i^T S(q_i)\gamma_i$, and $q_i^T\omega_i=0$, the variation of the action integral given at \refeqn{delG3} can be rewritten as
\begin{align}
\delta \tilde{\mathfrak{G}} & = 
\sum_{i=1}^n\int_{t_0}^{t_f} \Big\{ -\dot\pi_i+S(\omega_i)\pi_i -S(q_i)\deriv{\tilde H(q,\pi)}{q_i}\\
&\quad +S(q_i)\pi_i q_i^T\deriv{\tilde H(q,\pi)}{\pi_i}\Big\}\cdot  \gamma_i\nonumber\\
&\quad+\braces{\omega_i-(I_{3\times 3}-q_iq_i^T)\deriv{\tilde H(q,\pi)}{\pi_i} }\cdot\delta\pi_i^M\, dt.
\label{eqn:delG4}
\end{align}
According to Hamilton's phase space variational principle, $\delta\mathfrak{G}=-\delta\mathcal{W}$, the expression at the first pair of braces of the above equation, added with $S(q_i)f_i$, should be parallel to $q_i$, or equivalently,
\begin{gather}
(I_{3\times 3}-q_iq_i^T) \Big\{-\dot\pi_i+S(\omega_i)\pi_i -S(q_i)\deriv{\tilde H(q,\pi)}{q_i}\nonumber\\+S(q_i)\pi_i q_i^T\deriv{\tilde H(q,\pi)}{\pi_i}+S(q_i)f_i\Big\}=0.\label{eqn:delG4a}
\end{gather}
Also, the expression at the second pair of braces of \refeqn{delG4} is already parallel to $q_i$ to yield
\begin{gather*}
\omega_i = (I_{3\times 3}-q_iq_i^T)\deriv{\tilde H(q,\pi)}{\pi_i}=-S(q_i)^2\deriv{\tilde H(q,\pi)}{\pi_i}.
\end{gather*}
Substituting the second equation into \refeqn{rotkin5} yields \refeqn{5HEpifa}. 

Using the facts that $(I_{3\times 3}-q_iq_i^T) S(q_i) = -S(q_i)^3 =S(q_i)$ and $(I_{3\times 3}-q_iq_i^T) S(\omega_i)\pi_i = -S(q_i)\{ q_i\times (\omega_i\times \pi_i)\} = -S(q_i)\{ (q_i\cdot\pi_i)\omega_i - (q_i\cdot\omega_i)\pi_i\} = 0$, \refeqn{delG4a} reduces to
\begin{align*}
-(I_{3\times 3}-q_iq_i^T)\dot\pi_i -S(q_i)\parenth{\deriv{\tilde H}{q_i}-f_i}+S(q_i)\pi_i q_i^T\deriv{\tilde H}{\pi_i}=0.
\end{align*}
However, this is incomplete since it only determines the component of $\dot\pi_i$ that is perpendicular to $q_i$. The component of $\dot\pi_i$ that is parallel to $q_i$ is determined by taking the time derivative of $q_i\cdot\pi_i=0$ to obtain $q_i\cdot \dot\pi_i =-\dot q_i\cdot\pi_i$. Therefore, $q_iq_i^T\dot\pi_i = -q_i\pi_i^T\dot q_i$. By combining these, 
\begin{align*}
\dot\pi_i & =  -S(q_i)\parenth{\deriv{\tilde H}{q_i}-f_i}+S(q_i)\pi_i q_i^T\deriv{\tilde H}{\pi_i}+q_i\pi_i^TS(q_i)\deriv{\tilde H}{\pi_i}\nonumber\\
& = -S(q_i)\parenth{\deriv{\tilde H}{q_i}-f_i}+\deriv{\tilde H}{\pi_i}\times ((S(q_i)\pi_i)\times q_i)\nonumber\\
& = -S(q_i)\parenth{\deriv{\tilde H}{q_i}-f_i}+\deriv{\tilde H}{\pi_i}\times (-S(q_i)^2\pi_i).
\end{align*}
But, $-S(q_i)^2\pi_i=\pi_i$ since $\pi_i$ is normal to $q_i$. This yields \refeqn{5HEpifb}.

Next, \refeqn{5HEpigg} is obtained by substituting the angular velocity \refeqn{wi3} into the kinematics equation \refeqn{rotkin5}, and \refeqn{5HEpihh} is derived by substituting the Hamiltonian \refeqn{5Hpi} into \refeqn{5HEpifb}.

\bibliography{BibMaster}
\bibliographystyle{IEEEtran}

\end{document}